\documentclass[11pt]{amsart}
\usepackage{setspace}
\usepackage{mathrsfs}
\newtheorem{thm}{Theorem}[section]

\theoremstyle{definition}
\newtheorem{defn}[thm]{Definition}
\theoremstyle{remark}
\newtheorem{rem}[thm]{Remark}

\begin{document}

\title[Time-inhomogeneous random walks]{Conditions for recurrence and transience for time-inhomogeneous random walks} %
\author[Abramov]{Vyacheslav M. Abramov}
\address{24 Sagan Drive, Cranbourne North, Victoria-3977, Australia}%

\email{vabramov126@gmail.com}%


\begin{abstract}
The present paper extends the earlier results obtained in [V. M. Abramov, Bull. Aust. Math. Soc. \textbf{109} (2024), 393--402] for the case of time-inhomogeneous random walks, the increments of which take values in $\mathbb{R}$. By this, we give a full solution of the open problem formulated in [M. Menshikov and S. Volkov, Electron. J. Probab., \textbf{13} (2008), paper No. 31, 944--960] that was partially solved in the aforementioned paper by Abramov.
\end{abstract}

\maketitle

\section{Introduction}\label{S0}

The present paper is stimulated by Menshikov and Volkov \cite{MV}, where the authors studied the time-inhomogeneous random walk $X_n\in\mathbb{R}$, $n=1,2,\ldots$, with the drift satisfying the condition
\begin{equation}\label{1}
\mathsf{E}\{X_{n+1}-X_n~|~X_n=x\}\sim \rho \frac{|x|^\alpha}{n^\beta},
\end{equation}
with some constants $\rho>0$, $\alpha$ and $\beta$ and the special meaning of the symbol `$\sim$' further clarified for the case studies of the paper. The study in \cite{MV} was conducted under special assumptions, and one of the case studies, when $\alpha=2\beta-1$ and $0<\beta<1$ was left undecided. This case was partially resolved in recent paper \cite{A}, where process $X_t$ took the values in $\mathbb{Z}_+$ and fell into the category of time-inhomogeneous birth-and-death processes.

The aim of the present paper is to extend the results of \cite{A} for the case of the time-inhomogeneous random walk $X_n$, the increments of which belong to $\mathbb{R}$ and thus fully resolve the open problem from \cite{MV}.

It was believed in \cite{A} that the reduction from the case of a random walk taking the values in $\mathbb{Z}_+$ to the case of the random walk taking the values in $\mathbb{R}_+$ can be derived by a standard way. Further discussion of this problem with Dr. Tuan-Minh Nguyen, convinced us that the aforementioned reduction is not elementary, and therefore the solution to the open problem given in \cite{A} cannot be considered as a full solution. Indeed, the direct way of the proof anticipates further extension of \cite[Lemma 3.1]{A}, the way of which is unclear. As well, it turns out from the study of the present paper that such extension involves the theory of point processes and questions of measurability for the martingales that arise in the solution of the problem.
This motivated our present study.

In this study, we derived the system of equations for limit distribution closely related to that had been appeared in our previous study in \cite{A} for time-inhomogeneous birth-and-death processes. As well, the study is based on new ideas that recently have been developed in the paper on Markov chains \cite{A1}.
With the aid of these studies, we reduce our results to the earlier known ones from \cite{A}. For past studies of recurrence, transience or stability of time-inhomogeneous Markov processes see, e.g., \cite{AL, Z, ZI}).

Our approach is fully based on stochastic calculus and point processes theory, the means of which helped us to derive the stochastic variant of Chapman-Kolmogorov equations and provide their further analysis by the means of the theory of martingales in continuous time. The general method of the derivation of the stochastic variant of Chapman-Kolmogorov equations goes back to the paper by Kogan and Liptser \cite{KL}, where it was used for the first time.
\smallskip

The rest of the paper is organized as follows. The main result is formulated in Section \ref{S2}. The proof is given in Section \ref{S4}. We conclude the paper in Section \ref{S6}.

\section{The main result}\label{S2}
\subsection{The random walk}\label{S2.1}
We consider the continuous time process $Z(\tau)=\Lambda(\tau)-\mathrm{M}(\tau)$ taking values in $\mathbb{R}$, where the processes $\Lambda(\tau)$ and $\mathrm{M}(\tau)$ are time-inhomogeneous compound Poisson processes continuously depending on state and time, $\Lambda(0)=\mathrm{M}(0)=0$. The instantaneous rates of these processes are $\lambda_{Z(\tau),\tau}$ and $\mu_{Z(\tau),\tau}$, respectively, $(\lambda_{Z(\tau),\tau}+\mu_{Z(\tau),\tau}=1)$. The assumption $\lambda_{Z(\tau),\tau}+\mu_{Z(\tau),\tau}=1$ will enable us to consider the indexed parameters $\lambda$ and $\mu$ as conditional probabilities, provided that $\tau$ and $Z(\tau)$ are given.

The jumps of the process $\Lambda(\tau)$ are assumed to be positive independent identically distributed (i.i.d.) random variables $X_1$, $X_2$,\ldots, with mean 1 and finite variance $\mathsf{Var}(X_1)$.
The jumps of the process $\mathrm{M}(\tau)$ are assumed to be i.i.d. positive random variables $Y_1$, $Y_2$,\ldots, with mean 1 and finite variance $\mathsf{Var}(Y_1)$.


As well, the processes $\Lambda(\tau)$ and $\mathrm{M}(\tau)$ are assumed to be independent.
In the case when $X_1=Y_1=1$, the process reduces to the time-inhomogeneous birth-and-death process studied in \cite{A}.

We assume
\begin{eqnarray}
\lim_{\tau\to\infty}\lambda_{x,\tau}&=&\lambda_x>0,\label{3}\\
\lim_{\tau\to\infty}\mu_{x,\tau}&=&\mu_x>0,\label{5}
\end{eqnarray}
for all $x\in\mathbb{R}$.

\smallskip
The Markov process $Z(\tau)$ is a continuous time Markov process with continuous set of states. For recurrence of this process we have the following definition.

\begin{defn}
We call the Markov process $Z(\tau)$ \textit{recurrent}, if for any nonempty open set $\mathscr{A}\subset\mathbb{R}$ there exists an infinite sequence of times $\tau_1$, $\tau_2$,\ldots such that $Z(\tau_i)\in\mathscr{A}$, $i=1,2,\ldots$.
\end{defn}

The random walk that have been described here is a continuous time stochastic process, and the main result of this paper will be provided in terms of this process. Possible reduction to the random walk in discrete time, if we study the process in the points of jumps only, can be found in \cite{A} and in the relevant references provided there.

\subsection{Formulation of the main result}\label{S2.2}
The main result of this paper is the theorem formulated below.

Let $\tau_1$, $\tau_2$,\ldots denote the times of consecutive jumps of the random walk. We consider the random walk given by
\begin{eqnarray}
\mathsf{P}\{Z(\tau_{i+1})=x+X_i~|~Z(\tau_i)=x\}&=&\frac{1}{2}+\varphi(x,\tau_i),\label{4}\\
\mathsf{P}\{Z(\tau_{i+1})=x-Y_i~|~Z(\tau_i)=x\}&=&\frac{1}{2}-\varphi(x,\tau_i),\label{6}
\end{eqnarray}
where the function $\varphi(x,\tau_i)$ well-defines the right-hand sides of
\eqref{4} and \eqref{6} (i.e. keeps the probability distributions correctly defined).

The function $\varphi(x, \tau_i)$ is assumed to satisfy the following properties:
\begin{enumerate}
\item[(i)] the function $\varphi(x,\tau)$ is a positive Borel function,
\item[(ii)] $\varphi(\mathsf{E}Z(\tau), \tau)$ is decreasing in $\tau$, and $\varphi(x,\tau)\to0$  as $x\to\infty$,
\item[(iii)] for any $\epsilon>0$ there exist $x_0$ and $\delta>0$ such that for all $x>x_0$
\begin{eqnarray*}
\big|\varphi\big(x, x^2(1+\delta)\big)-\varphi(x, x^2)\big|&<&\epsilon\varphi(x, x^2),\\
\big|\varphi\big(x, x^2(1-\delta)\big)-\varphi(x, x^2)\big|&<&\epsilon\varphi(x, x^2).
\end{eqnarray*}
\end{enumerate}

\begin{thm}\label{T1}
Let  $\lambda_{x,t}=1/2+\varphi(x,t)$, $\mu_{x,t}=1/2-\varphi(x,t)$. Then, the random walk is recurrent, if there exists $x_0>0$ and $c<1$ such that for all $x\geq x_0$  we have $\varphi(x,x^2)\leq c/(4x)$, and it is transient if there exists $x_0>0$ and $c>1$ such that for all $x\geq x_0$ we have $\varphi(x,x^2)\geq c/(4x)$.
\end{thm}

\begin{rem}
Following \eqref{1}, Menshikov and Volkov \cite{MV} describe the asymptotic behavior of the random walk via the conditional expectations of its increments. In the present study, we describe the asymptotic behavior of the random walk via conditional probabilities written in terms of intensities (see relations \eqref{4} and \eqref{6}). As well, the basic assumptions, under which the study is conducted in \cite{MV} and in the present paper are different. So, both of the studies have been provided in different terms, and the solution to the open problem of Menshikov and Volkov \cite{MV} is given here in the terms of the present study.
For instance, in the case where
\[
\varphi(x,n)=\lambda_{x,n}-\frac{1}{2}=\frac{1}{2}-\mu_{x,n}=\frac{1}{2}\rho|x|^\alpha n^{-\beta},
\]
the process $Z(\tau)$ is more particular compared to the process considered by Menshikov and Volkov \cite{MV}, since the Menshikov--Volkov random walk does not require that increments are generated by identically distributed random variables. Unfortunately we cannot relax our assumptions on the sequences of random variables $\{X_i\}_{i\geq1}$ and $\{Y_i\}_{i\geq1}$ in the direct way.  However, possible extensions of the considered random walk that enables us to obtain our result under weaker assumptions can be provided via suitable approximations, and it will be discussed later in the paper (see Section \ref{S3.3}).
\end{rem}

\section{Proof of the main result}\label{S4}
\subsection{Description of discrete valued processes}\label{S4.1}
In our study, the continuous time process $Z(\tau)=\Lambda(\tau)-\mathrm{M}(\tau)$ taking values in $\mathbb{R}$ will be approximated (in the sense of the pointwise convergence of discrete distributions to their continuous limits indicated exactly in Section \ref{S4.3})  by the series of processes $Z^{(q)}(\tau)=\Lambda^{(q)}(\tau)-\mathrm{M}^{(q)}(\tau)$, where the meaning of the processes $\Lambda^{(q)}(\tau)$ and $\mathrm{M}^{(q)}(\tau)$ is similar to the corresponding processes $\Lambda(\tau)$ and $\mathrm{M}(\tau)$, with the only difference that the jumps $X_1^{(q)}$, $X_2^{(q)}$,\ldots and $Y_1^{(q)}$, $Y_2^{(q)}$,\ldots take discrete set of values. Specifically, we assume that the jumps $X_1^{(q)}$ and $Y_1^{(q)}$ take the values $jq$, $j=1,2,\ldots$, keeping the properties $\mathsf{E}X_1^{(q)}=\mathsf{E}Y_1^{(q)}=1$ and $\mathsf{Var}(X_1^{(q)})<\infty$, $\mathsf{Var}(Y_1^{(q)})<\infty$.  These properties are assumed to satisfy for all $q\leq q_0$, and for vanishing $q$ we have the properties
\[\lim_{q\downarrow0}\mathsf{E}X_1^{(q)}=\lim_{q\downarrow0}\mathsf{E}Y_1^{(q)}=1,\] and
\[\lim_{q\downarrow0}\mathsf{Var}(X_1^{(q)})=\mathsf{Var}(X_1)<\infty, \quad \lim_{q\downarrow0}\mathsf{Var}(Y_1^{(q)})=\mathsf{Var}(Y_1)<\infty.\]

Below we define the processes $\Lambda_{nq,\tau}^{(q)}(\tau)$ and $\mathrm{M}_{nq,\tau}^{(q)}(\tau)$ using the scheme given in \cite{A1}. The sub-indices $nq$ and $\tau$ emphasise the dependence on state $nq$ and time $\tau$. Note also that with given $q$,  the process $Z^{(p)}(\tau)$ takes the values on the infinite lattice $\ldots$, $-2q$, $-q$, $0$, $q$, $2q$, $\ldots$.

\subsubsection*{Definition of the processes $\Lambda_{nq,\tau}^{(q)}(\tau)$ and $\mathrm{M}_{nq,\tau}^{(q)}(\tau)$}
Let $\tau_1$, $\tau_2$,\ldots be the times of jumps, when $\mathsf{I}\{Z^{(q)}(\tau_i-)=nq, Z^{(q)}(\tau_i)>nq\}=1$. These times of jumps are denoted by $\tau_1^{(nq)}$, $\tau_2^{(nq)}$,\ldots.
Similarly, let $\eta_1$, $\eta_2$,\ldots be the times of jumps, when $\mathsf{I}\{Z^{(q)}(\eta_i-)=nq, Z^{(q)}(\eta_i)<nq\}=1$. These times of jumps are denoted by $\eta_1^{(nq)}$, $\eta_2^{(nq)}$,\ldots. Denote by $\triangle Z^{(q)}(\tau_i^{(nq)})=Z^{(q)}(\tau_i^{(nq)})-Z^{(q)}(\tau_i^{(nq)}-)$ the jumps at times $\tau_i^{(nq)}$, $i=1,2,\ldots$, and by $\triangle Z^{(q)}(\eta_i^{(nq)})=Z^{(q)}(\eta_i^{(nq)})-Z^{(q)}(\eta_i^{(nq)}-)$ the jumps at times $\eta_i^{(nq)}$, $i=1,2,\ldots$.
The sequence $\triangle Z^{(q)}(\tau_i^{(nq)})$ (containing positive values) is denoted by $\triangle\Lambda_{nq,\tau}^{(q)}(\tau_i^{(nq)})$ and associated with the process $\Lambda_{nq,\tau}^{(q)}(\tau)$, where the process $\Lambda_{nq,\tau}^{(q)}(\tau)$ itself is obtained by deleting all the intervals of the process $\Lambda^{(q)}(\tau)$ where $\mathsf{I}\{Z^{(q)}(\tau-)\neq nq\}$ and merging the ends, with following re-scaling times. The sequence $\triangle Z^{(q)}(\eta_i^{(nq)})$ (containing negative values) taken with opposite (positive) sign is denoted by $\triangle\mathrm{M}_{nq,\tau}^{(q)}(\eta_i^{(nq)})$ and associated with the process $\mathrm{M}_{nq,\tau}^{(q)}(\tau)$, where the process $\mathrm{M}_{nq,\tau}^{(q)}(\tau)$ is obtained similarly to that of $\Lambda_{nq,\tau}^{(q)}(\tau)$. The elements of the sequences $\triangle\Lambda_{nq,\tau}^{(q)}(\tau_i^{(nq)})$ and $\triangle\mathrm{M}_{nq,\tau}^{(q)}(\eta_i^{(nq)})$ are renumbered and denoted by $X_1^{(n,q)}$, $X_2^{(n,q)}$,\ldots and $Y_1^{(n,q)}$, $Y_2^{(n,q)}$,\ldots, respectively. The sequences $\{X_i^{(n,q)}\}_{i\geq1}$ and $\{Y_i^{(n,q)}\}_{i\geq1}$ are independent and consist of i.i.d. random variables; $X_1^{(n,q)}$ has the same distribution as $X_1^{(q)}$, $Y_1^{(n,q)}$ has the same distribution as $Y_1^{(q)}$, and the double upper indexation for $X_i^{(n,q)}$ and $Y_i^{(n,q)}$ is merely used for correct numbering of the random variables in the sequences in order to avoid using same random variables for different $n$.

For $\Lambda_{nq,\tau}^{(q)}(\tau)$ and $\mathrm{M}_{nq,\tau}^{(q)}(\tau)(\tau)$ we have:
\[
\Lambda_{nq,\tau}^{(q)}(\tau)=\sum_{i=1}^{N_{nq,\tau}^{\Lambda}(\tau)}X_i^{(n,q)}=\int_{0}^{\tau}X_{N_{nq,t}^\Lambda(t)}^{(n,q)}\mathrm{d}N_{nq,t}^\Lambda(t)
\]
and
\[
\mathrm{M}_{nq,\tau}^{(q)}(\tau)=\sum_{i=1}^{N_{nq,\tau}^{\mathrm{M}}(\tau)}Y_i^{(n,q)}=\int_{0}^{\tau}Y_{N_{nq,t}^\mathrm{M}(t)}^{(n,q)}\mathrm{d}N_{nq,t}^\mathrm{M}(t),
\]
where $N_{nq,\tau}^{\Lambda}(\tau)$ and $N_{nq,\tau}^{\mathrm{M}}(\tau)$ are underlying point process of jump times.

 As $q$ vanishes (remaining positive) and $n\to\pm\infty$, we have the required approximation for the limiting distribution of $Z(\tau)$ as $\tau$ increases to infinity. (Here and further $n\to\pm\infty$ means that either $n\to\infty$ or $n\to-\infty$. That is, for a fixed $x$, $nq\to x$ implies $q\downarrow0$, $n\to+\infty$ if $x$ positive, or $q\downarrow0$, $n\to-\infty$ if $x$ negative.) Note also
 that the instantaneous rates of the Poisson processes $\Lambda_{jq,\tau}^{(q)}(\tau)$ and $\mathrm{M}_{jq,\tau}^{(q)}(\tau)$ are denoted by $\lambda_{jq,\tau}^{(q)}$ and $\mu_{jq,\tau}^{(q)}$, respectively ($\lambda_{jq,\tau}^{(q)}+\mu_{jq,\tau}^{(q)}=1$). Assumptions \eqref{3} and \eqref{5} in this case are rewritten
\begin{eqnarray*}
\lim_{\tau\to\infty}\lambda_{nq,\tau}^{(q)}&=&\lambda_{nq}^{(q)}>0,\\
\lim_{\tau\to\infty}\mu_{nq,\tau}^{(q)}&=&\mu_{nq}^{(q)}>0.
\end{eqnarray*}

\subsection{Stochastic equations and semimartingale decompositions}\label{S4.2}
The system of equations for the indicators $\mathsf{I}\{Z^{(q)}(\tau)=nq\}$, $n=...,-1,0,1,\ldots$, that holds with probability $1$ is
\begin{equation}\label{7}
\begin{aligned}
\mathsf{I}\{Z^{(q)}(\tau)=nq\}&=\sum_{j=n+1}^{\infty}\mathsf{I}\{Z^{(q)}(\tau-)=jq\}\mathsf{I}\{\triangle\mathrm{M}^{(q)}_{jq,\tau}(\tau)=(j-n)q\}\\
+&\sum_{j=-\infty}^{n-1}\mathsf{I}\{Z^{(q)}(\tau-)=jq\}\mathsf{I}\{\triangle\Lambda_{jq,\tau}^{(q)}(\tau)=(n-j)q\}\\
+&\mathsf{I}\{Z^{(q)}(\tau-)=nq\}\mathsf{I}\{\triangle\mathrm{M}_{nq,\tau}^{(q)}(\tau)=0\}\mathsf{I}\{\triangle\Lambda_{nq,\tau}^{(q)}(\tau)=0\}.
\end{aligned}
\end{equation}

To derive the expressions for the differentials of indicator functions that involve the processes $\Lambda_{jq,\tau}^{(q)}(\tau)$ and $\mathrm{M}_{jq,\tau}^{(q)}(\tau)$ we need to express these processes via the point processes $N_{jq,\tau}^{\Lambda}(\tau)$ and $N_{jq,\tau}^{\mathrm{M}}(\tau)$.

Note that
\begin{eqnarray*}
\triangle \Lambda_{jq,\tau}^{(q)}(\tau)&=&X_{N_{jq,\tau}^\Lambda(\tau)}^{(j,q)}\triangle N_{jq,\tau}^\Lambda(\tau),\\
\triangle \mathrm{M}_{jq,\tau}^{(q)}(\tau)&=&Y_{N_{jq,\tau}^\mathrm{M}(\tau)}^{(j,q)}\triangle N_{jq,\tau}^\mathrm{M}(\tau),
\end{eqnarray*}
where $\triangle N_{jq,\tau}^\Lambda(\tau)$ and $\triangle N_{jq,\tau}^\mathrm{M}(\tau)$ denote the increments of the corresponding point processes at time $\tau$ taking the values $0$ (no jump at time $\tau$ occurs) or $1$ (a jump occurs exactly at time $\tau$).

We have the following relations:
\begin{align}
  \mathsf{I}\{\triangle\Lambda_{jq,\tau}^{(q)}(\tau)=(n-j)q\} &=\mathsf{I}\{X_{N_{jq,\tau}^\Lambda(\tau)}^{(j,q)}=(n-j)q\}\triangle N_{jq,\tau}^\Lambda(\tau)\label{30}  \\
   & =\mathsf{I}\{X_{N_{jq,\tau}^\Lambda(\tau-)+1}^{(j,q)}=(n-j)q\}\triangle N_{jq,\tau}^\Lambda(\tau),\label{31}
\end{align}
and
\begin{align}
  \mathsf{I}\{\triangle\mathrm{M}_{jq,\tau}^{(q)}(\tau)=(j-n)q\} &=\mathsf{I}\{Y_{N_{jq,\tau}^\mathrm{M}(\tau)}^{(j,q)}=(j-n)q\}\triangle N_{jq,\tau}^\mathrm{M}(\tau)\label{32}  \\
   & =\mathsf{I}\{Y_{N_{jq,\tau}^\mathrm{M}(\tau-)+1}^{(j,q)}=(j-n)q\}\triangle N_{jq,\tau}^\mathrm{M}(\tau).\label{33}
\end{align}
Representations \eqref{30} and \eqref{31} are equivalent, since
\[
\mathsf{I}\{X_{N_{jq,\tau}^\Lambda(\tau)}^{(j,q)}=(n-j)q\}\triangle N_{jq,\tau}^\Lambda(\tau)=\mathsf{I}\{X_{N_{jq,\tau}^\Lambda(\tau-)+1}^{(j,q)}=(n-j)q\}\triangle N_{jq,\tau}^\Lambda(\tau)
\]
if $\triangle N_{jq,\tau}^\Lambda(\tau)=1$. The same observation pertains to \eqref{32} and \eqref{33}. The significance of these relations will be seen later when we discuss the predictability of $\mathsf{I}\{X_{N_{jq,\tau}^\Lambda(\tau-)+1}^{(j,q)}=(n-j)q\}$ and $\mathsf{I}\{Y_{N_{jq,\tau}^\mathrm{M}(\tau-)+1}^{(j,q)}=(j-n)q\}$ for the martingale arguments used in the paper. (For the further measurability arguments, we need to specify all the processes that appear in the stochastic equations at time $\tau-$, that is, immediately before a possible jump at time $\tau$.)

In addition to these equation, with probability $1$ we have
\begin{equation}\label{34}
\begin{aligned}
&\mathsf{I}\{\triangle\Lambda_{nq,\tau}^{(q)}(\tau)=0\}\mathsf{I}\{\triangle\mathrm{M}_{nq,\tau}^{(q)}(\tau)=0\}\\
&=[1-X_{N_{nq,\tau}^\Lambda(\tau)}^{(n,q)}\triangle N_{nq,\tau}^\Lambda(\tau)][1-Y_{N_{nq,\tau}^\mathrm{M}(\tau)}^{(n,q)}\triangle N_{nq,\tau}^\mathrm{M}(\tau)]\\
&=[1-X_{N_{nq,\tau}^\Lambda(\tau-)}^{(n,q)}\triangle N_{nq,\tau}^\Lambda(\tau)][1-Y_{N_{nq,\tau}^\mathrm{M}(\tau-)}^{(n,q)}\triangle N_{nq,\tau}^\mathrm{M}(\tau)].
\end{aligned}
\end{equation}
Here in \eqref{34} we use the fact that $\tau$ is not the point of jump, and hence $N_{nq,\tau}^\Lambda(\tau)=N_{nq,\tau}^\Lambda(\tau-)$, $N_{nq,\tau}^\mathrm{M}(\tau)=N_{nq,\tau}^\mathrm{M}(\tau-)$ that results in the corresponding relations $X_{N_{nq,\tau}^\Lambda(\tau)}^{(n,q)}=X_{N_{nq,\tau}^\Lambda(\tau-)}^{(n,q)}$, $Y_{N_{nq,\tau}^\mathrm{M}(\tau)}^{(n,q)}=Y_{N_{nq,\tau}^\mathrm{M}(\tau-)}^{(n,q)}$ used above.

To derive stochastic differential equations, relations \eqref{31}, \eqref{33} and \eqref{34} are rewritten below in the terms of stochastic differentials. From \eqref{31} and \eqref{33} we obtain:
\begin{equation}\label{35}
\mathrm{d}\mathsf{I}\{\triangle\Lambda_{jq,\tau}^{(q)}(\tau)=(n-j)q\}=\mathsf{I}\{X_{N_{jq,\tau}^\Lambda(\tau-)+1}^{(j,q)}=(n-j)q\}\mathrm{d} N_{jq,\tau}^\Lambda(\tau),
\end{equation}

\begin{equation}\label{36}
\mathrm{d}\mathsf{I}\{\triangle\mathrm{M}_{jq,\tau}^{(q)}(\tau)=(j-n)q\}=\mathsf{I}\{Y_{N_{jq,\tau}^\mathrm{M}(\tau-)+1}^{(j,q)}=(j-n)q\}\mathrm{d} N_{jq,\tau}^\mathrm{M}(\tau),
\end{equation}

In turn, from \eqref{34} we obtain
\begin{equation}\label{37}
\begin{aligned}
&\mathrm{d}(\mathsf{I}\{\triangle\Lambda_{nq,\tau}^{(q)}(\tau)=0\}\mathsf{I}\{\triangle\mathrm{M}_{nq,\tau}^{(q)}(\tau)=0\})\\
&=-X_{N_{nq,\tau}^\Lambda(\tau-)}^{(n,q)}\mathrm{d}N_{nq,\tau}^\Lambda(\tau)-Y_{N_{nq,\tau}^\mathrm{M}(\tau-)}^{(n,q)}\mathrm{d}N_{nq,\tau}^\mathrm{M}(\tau).
\end{aligned}
\end{equation}

According to representations \eqref{35}--\eqref{37}, the differential form of equation \eqref{7} is
\begin{equation}\label{8}
\begin{aligned}
&\mathrm{d}\mathsf{I}\{Z^{(q)}(\tau)=nq\}\\
=&\sum_{j=n+1}^{\infty}\mathsf{I}\{Z^{(q)}(\tau-)=jq\}\mathsf{I}\{Y_{N_{jq,\tau}^{\mathrm{M}}(\tau-)+1}^{(j,q)}=(j-n)q\}\mathrm{d}N_{jq,\tau}^{\mathrm{M}}(\tau)\\
&+\sum_{j=-\infty}^{n-1}\mathsf{I}\{Z^{(q)}(\tau-)=jq\}\mathsf{I}\{X_{N_{jq,\tau}^{\Lambda}(\tau-)+1}^{(j,q)}=(n-j)q\}\mathrm{d}N_{jq,\tau}^{\Lambda}(\tau)\\
&-\mathsf{I}\{Z^{(q)}(\tau-)=nq\}[X_{N_{nq,\tau}^\Lambda(\tau-)}^{(n,q)}\mathrm{d}N_{nq,\tau}^\Lambda(\tau)+Y_{N_{nq,\tau}^\mathrm{M}(\tau-)}^{(n,q)}\mathrm{d}N_{nq,\tau}^\mathrm{M}(\tau)],
\end{aligned}
\end{equation}
with probability $1$.

In the sequel, we need in the Doob-Meyer semimartingale decomposition (see \cite{J-Sh, L-Sh, Pr}) for the processes $N^\Lambda_{jq,\tau}(\tau)$ and $N^\mathrm{M}_{jq,\tau}(\tau)$ with respect to the filtration $\mathcal{F}_\tau^{(q)}$ that is defined as
\[
\mathcal{F}_\tau^{(q)}=\sigma\{Z^{(q)}(t), 0\leq t\leq\tau; X_{i}^{(q,j)}, Y_k^{(q,j)}, 1\leq i,k<\infty, -\infty<j<\infty \}.
\]
For this filtration, the sample paths of $Z^{(q)}(t)$, $N_{jq,t}^\Lambda$, $N_{jq,t}^\mathrm{M}$, $-\infty<j<\infty$ are known up to time $\tau$, and all future jump sizes after time $\tau$ of $N_{jq,t}^\Lambda$, $N_{jq,t}^\mathrm{M}$, $-\infty<j<\infty$ are known, but the locations of jumps are unknown.

Since the instantaneous rates of the processes $N^\Lambda_{jq,\tau}(\tau)$ and $N^\mathrm{M}_{jq,\tau}(\tau)$ at time $\tau$ coincide with the instantaneous rates of the ordinary Poisson processes having the parameters $\lambda_{jq,\tau}^{(q)}$ and $\mu_{jq,\tau}^{(q)}$, respectively, then the Doob-Meyer semimartingale decompositions can be written in the form of stochastic differentials:
\begin{eqnarray}
\mathrm{d}N^\Lambda_{jq,\tau}(\tau)&=&\lambda_{jq,\tau}^{(q)}\mathrm{d}\tau+\mathrm{d}\boldsymbol{M}_{N^\Lambda_{jq,\tau}(\tau)},\label{9}\\
\mathrm{d}N^\mathrm{M}_{jq,\tau}(\tau)&=&\mu_{jq,\tau}^{(q)}\mathrm{d}\tau+\mathrm{d}\boldsymbol{M}_{N^\mathrm{M}_{jq,\tau}(\tau)},\label{10}
\end{eqnarray}
where $\boldsymbol{M}_{N^\Lambda_{jq,\tau}(\tau)}$ and $\boldsymbol{M}_{N^\mathrm{M}_{jq,\tau}(\tau)}$ denote square integrable martingales in the decompositions, $j=\ldots,-2,-1,0,1,2,\ldots$

 Then the system of equations \eqref{8} is rewritten as follows:
\begin{equation}\label{15}
\begin{aligned}
&\mathrm{d}\mathsf{I}\{Z^{(q)}(\tau)=nq\}\\
=&\sum_{j=n+1}^{\infty}\mathsf{I}\{Z^{(q)}(\tau-)=jq\}\mathsf{I}\{Y_{N_{jq,\tau}^{\mathrm{M}}(\tau-)+1}^{(j,q)}=(j-n)q\}\mu_{jq,\tau}^{(q)}\mathrm{d}\tau\\
&+\sum_{j=-\infty}^{n-1}\mathsf{I}\{Z^{(q)}(\tau-)=jq\}\mathsf{I}\{X_{N_{jq,\tau}^{\Lambda}(\tau-)+1}^{(j,q)}=(n-j)q\}\lambda_{jq,\tau}^{(q)}\mathrm{d}\tau\\
&-\mathsf{I}\{Z^{(q)}(\tau-)=nq\}[\lambda_{nq,\tau}^{(q)}+\mu_{nq,\tau}^{(q)}]\mathrm{d}\tau\\
&+\sum_{j=n+1}^{\infty}\mathsf{I}\{Z^{(q)}(\tau-)=jq\}\mathsf{I}\{Y_{N_{jq,\tau}^{\mathrm{M}}(\tau-)+1}^{(j,q)}=(j-n)q\}\mathrm{d}\boldsymbol{M}_{N^\mathrm{M}_{jq,\tau}(\tau)}\\
&+\sum_{j=-\infty}^{n-1}\mathsf{I}\{Z^{(q)}(\tau-)=jq\}\mathsf{I}\{X_{N_{jq,\tau}^{\Lambda}(\tau-)+1}^{(j,q)}=(n-j)q\}\mathrm{d}\boldsymbol{M}_{N^\Lambda_{jq,\tau}(\tau)}\\
&-\mathsf{I}\{Z^{(q)}(\tau-)=nq\}[X_{N_{nq,\tau}^{\Lambda}(\tau-)}^{(n,q)}\mathrm{d}\boldsymbol{M}_{N^\Lambda_{nq,\tau}(\tau)}+Y_{N_{nq,\tau}^{\mathrm{M}}(\tau-)}^{(n,q)}\mathrm{d}\boldsymbol{M}_{N^\mathrm{M}_{nq,\tau}(\tau)}],
\end{aligned}
\end{equation}
where $\boldsymbol{M}_{N^\mathrm{M}_{jq,\tau}(\tau)}$ and $\boldsymbol{M}_{N^\Lambda_{jq,\tau}(\tau)}$, $-\infty<j<\infty$, all are $\mathcal{F}_\tau^{(q)}$-predictable.

 Now we rewrite the system of equations \eqref{15} in the integral form by taking the expectation and averaging. We have:

\begin{equation}\label{17}
\begin{aligned}
0=&\lim_{T\to\infty}\frac{1}{T}\mathsf{E}\int_0^T\sum_{j=n+1}^{\infty}\mathsf{I}\{Z^{(q)}(\tau-)=jq\}\mathsf{I}\{Y_{N_{jq,\tau}^{\mathrm{M}}(\tau-)+1}^{(j,q)}=(j-n)q\}\mu_{jq,\tau}^{(q)}\mathrm{d}\tau\\
&+\lim_{T\to\infty}\frac{1}{T}\mathsf{E}\int_0^T\sum_{j=-\infty}^{n-1}\mathsf{I}\{Z^{(q)}(\tau-)=jq\}\mathsf{I}\{X_{N_{jq,\tau}^{\Lambda}(\tau-)+1}^{(j,q)}=(n-j)q\}\lambda_{jq,\tau}^{(q)}\mathrm{d}\tau\\
&-\lim_{T\to\infty}\frac{1}{T}\mathsf{E}\int_0^T\mathsf{I}\{Z^{(q)}(\tau-)=nq\}[\lambda_{nq,\tau}^{(q)}+\mu_{nq,\tau}^{(q)}]\mathrm{d}\tau\\
&+\lim_{T\to\infty}\frac{1}{T}\mathsf{E}\int_0^T\sum_{j=n+1}^{\infty}\mathsf{I}\{Z^{(q)}(\tau-)=jq\}\mathsf{I}\{Y_{N_{jq,\tau}^{\mathrm{M}}(\tau-)+1}^{(j,q)}=(j-n)q\}\mathrm{d}\boldsymbol{M}_{N^\mathrm{M}_{jq,\tau}(\tau)}\\
&+\lim_{T\to\infty}\frac{1}{T}\mathsf{E}\int_0^T\sum_{j=-\infty}^{n-1}\mathsf{I}\{Z^{(q)}(\tau-)=jq\}\mathsf{I}\{X_{N_{jq,\tau}^{\Lambda}(\tau-)+1}^{(j,q)}=(n-j)q\}\mathrm{d}\boldsymbol{M}_{N^\Lambda_{jq,\tau}(\tau)}\\
&-\lim_{T\to\infty}\frac{1}{T}\mathsf{E}\int_0^T\mathsf{I}\{Z^{(q)}(\tau-)=nq\}[X_{N_{nq,\tau}^{\Lambda}(\tau-)}^{(n,q)}\mathrm{d}\boldsymbol{M}_{N^\Lambda_{nq,\tau}(\tau)}+Y_{N_{nq,\tau}^{\mathrm{M}}(\tau-)}^{(n,q)}\mathrm{d}\boldsymbol{M}_{N^\mathrm{M}_{nq,\tau}(\tau)}]\\
=&J_1+J_2+J_3+J_4+J_5+J_6,
\end{aligned}
\end{equation}
where all the limits in \eqref{17} exist. Indeed, the limits for $J_4$, $J_5$ and $J_6$ exist, since all them are predictable with respect to $\mathcal{F}^{(q)}=\cup_{\tau>0}\mathcal{F}_{\tau-}^{(q)}$ martingales with mean $0$, and hence $J_4=J_5=J_6=0$. The limits $J_1$, $J_2$, $J_3$ exist, since all the integrands are positive, and hence the integrals themselves are monotone increasing (cumulative) functions of the upper limit $T$.

For $J_1$ we have
 \begin{equation}\label{19}
 \begin{aligned}
J_1&=\lim_{T\to\infty}\frac{1}{T}\mathsf{E}\int_0^T\sum_{j=n+1}^{\infty}\mathsf{I}\{Z^{(q)}(\tau-)=jq\}\mathsf{I}\{Y_{N_{jq,\tau}^{\mathrm{M}}(\tau-)+1}^{(j,q)}=(j-n)q\}\mu_{jq,\tau}^{(q)}\mathrm{d}\tau\\
&=\lim_{T\to\infty}\frac{1}{T}\int_0^T\sum_{j=n+1}^{\infty}\mathsf{P}\{Z^{(q)}(\tau-)=jq\}p_Y^{(q)}((j-n)q)\mu_{jq,\tau}^{(q)}\mathrm{d}\tau\\
&=\mu_{nq+1}^{(q)}\lim_{T\to\infty}\mathsf{P}\{Z^{(q)}(T-)=nq+1\},
 \end{aligned}
 \end{equation}
where $p_Y^{(q)}(k)=\mathsf{P}\{Y_1^{(q)}=k\}$, $\mu_{nq}^{(q)}=\lim_{\tau\to\infty}\mu_{nq,\tau}^{(q)}$. Denoting also $p_X^{(q)}(k)=\mathsf{P}\{X_1^{(q)}=k\}$, $\lambda_{nq}^{(q)}=\lim_{\tau\to\infty}\lambda_{nq,\tau}^{(q)}$, for $J_2$ we have
\begin{equation}\label{20}
\begin{aligned}
J_2&=\lim_{T\to\infty}\frac{1}{T}\int_0^T\sum_{j=-\infty}^{n-1}\mathsf{P}\{Z^{(q)}(\tau-)=jq\}p_X^{(q)}((n-j)q)\lambda_{jq,\tau}^{(q)}\mathrm{d}\tau\\
&=\lambda_{nq-1}^{(q)}\lim_{T\to\infty}\mathsf{P}\{Z^{(q)}(T-)=nq-1\}.
\end{aligned}
\end{equation}
Similarly,
\begin{equation}\label{21}
\begin{aligned}
J_3&=-\lim_{T\to\infty}\frac{1}{T}\int_0^T\mathsf{P}\{Z^{(q)}(\tau-)=nq\}[\lambda_{nq,\tau}^{(q)}+\mu_{nq,\tau}^{(q)}]\mathrm{d}\tau\\
&=-[\lambda_{nq}^{(q)}+\mu_{nq}^{(q)}]\lim_{T\to\infty}\mathsf{P}\{Z^{(q)}(T-)=nq\}.
\end{aligned}
\end{equation}

\subsection{Discussion of the case when the sequences $\{X_i^{(q)}\}$ and $\{Y_i^{(q)}\}$ are state-dependent}\label{S3.3}
Our results are not changed, if we assume that the sequences $\{X_i^{(q)}\}$ and $\{Y_i^{(q)}\}$ are state-dependent rather than i.i.d., but having the same mean equal to $1$ and uniformly bounded variances.
In other words, we will assume that the sequences $\{X_i^{(n,q)}\}$ and $\{Y_i^{(n,q)}\}$ consist of i.i.d. random variables, however $X_1^{(n,q)}$ and $X_1^{(j,q)}$ as well as $Y_1^{(n,q)}$ and $Y_1^{(j,q)}$  have different distributions when $j\neq n$ in general.

In this case, the second line of \eqref{19} is slightly changed:
 \begin{equation*}\label{22}
 \begin{aligned}
J_1
&=\lim_{T\to\infty}\frac{1}{T}\int_0^T\sum_{j=n+1}^{\infty}\mathsf{P}\{Z^{(q)}(\tau-)=jq\}p_Y^{(j,q)}((j-n)q)\mu_{jq,\tau}^{(q)}\mathrm{d}\tau\\
 \end{aligned}
 \end{equation*}
 where $p_Y^{(j,q)}(k)=\mathsf{P}\{Y_1^{(j,q)}=k\}$. Since the expectations of $Y_1^{(j,q)}$ are equal to $1$ for all $j$, we finally arrive at the same result. The similar note is true for \eqref{20}. So, our results are insensitive to the form of distributions of $X_1^{(j,q)}$ and $Y_1^{(j,q)}$ for different $j$.




 \subsection{Connection with birth-and-death processes}\label{S4.3}
 From limits \eqref{19}--\eqref{21} we have the following system of equations:
 \begin{equation}\label{23}
 0=P_{nq+1}^{(q)}\mu_{nq+1}^{(q)}+P_{nq-1}^{(q)}\lambda_{nq-1}-P_{nq}^{(q)}(\lambda_{nq}^{(q)}+\mu_{nq}^{(q)}), \quad n\in\mathbb{Z},
 \end{equation}
 where $P_{nq}^{(q)}=\lim_{T\to\infty}\mathsf{P}\{Z^{(q)}(T-)=nq\}$. Now we need to take the limit as $nq\to x$. To do that let us study the properties of the functions $J_1$, $J_2$ and $J_3$ described in \eqref{19}--\eqref{21}. Note first that the summands in \eqref{19} and \eqref{20} are positive. That is, the convergence of the series implies their absolute convergence, and the orders of sum and integrals are interchangeable. In addition, for any $nq$ the integrands are bounded and uniformly continuous in $nq$, and the convergence to the limits is uniform. Taking limits in \eqref{23} as $nq\to x$ ($q\downarrow0$, $n\to\pm\infty$) implies interchange the order of the limits $nq\to x$ vs $T\to\infty$. Since the convergence in $T$ is uniform, then according to Moore-Osgood theorem \cite[pp. 139--140]{T} the required order of the limits can be interchanged. Then
 taking the limit as $nq\to x$ with $n\to\pm\infty$ and $q\downarrow0$ and integrating over integer intervals $(k-1,k]$, $k\in\mathbb{Z}$, we obtain:
 \begin{equation}\label{25}
 \begin{aligned}
 0&=\int_{k-1}^{k}P_{x+1}\mu_{x+1}\mathrm{d}x+\int_{k-1}^{k}P_{x-1}\lambda_{x-1}\mathrm{d}x-\int_{k-1}^{k}P_x(\lambda_x+\mu_x)\mathrm{d}x\\
 &=\mu_{k+1}^*\int_{k-1}^{n}P_{x+1}\mathrm{d}x+\lambda_{k-1}^*\int_{k-1}^{k}P_{x-1}\mathrm{d}x-(\lambda_k^*+\mu_k^*)\int_{k-1}^{n}P_x\mathrm{d}x.
 \end{aligned}
 \end{equation}
 In \eqref{25} we used the notation
 \[
 P_x=\lim_{\substack{n\to\pm\infty\\ q\downarrow0\\ nq\to x}}P_{nq}^{(q)}, \quad \lambda_x=\lim_{\substack{n\to\pm\infty\\ q\downarrow0\\ nq\to x}}\lambda_{nq}^{(q)}, \quad \mu_x=\lim_{\substack{n\to\pm\infty\\ q\downarrow0\\ nq\to x}}\mu_{nq}^{(q)}.
 \]

 The last relation in \eqref{25} is true due to the first mean value theorem for integration \cite[p. 159]{C}. Denoting $\int_{k-1}^{k}P_x\mathrm{d}x$ by $P_k^*$, we rewrite \eqref{25} in the form
 \begin{equation}\label{27}
 0=P_{k+1}^*\mu_{k+1}^*+P_{k-1}^*\lambda_{k-1}^*-P_k^*(\lambda_k^*+\mu_k^*), \quad k\in\mathbb{Z}.
 \end{equation}

 The infinite system of equations \eqref{27} describes the known system of equation for ordinary (time-homogeneous) birth-and-death process without boundary condition. That is, this system of equations describes a bilateral birth-and-death process \cite{P}. Under the setting $\lambda_k^*\geq\mu_k^*$, $k\in\mathbb{Z}$, this birth-and-death process is recurrent, if there exists $k_0>0$ such that $\lambda_k^*/\mu_k^*\leq1+k^{-1}$ for all $k\geq k_0$ and it is transient, if there exist $k_0>0$ and $c>1$ such that $\lambda_k^*/\mu_k^*\geq1+c/k$ for all $k\geq k_0$. This condition was used in \cite[Lemma 3.1]{A} for the ordinary birth-and-death process, but it remains true for bilateral as well, for which the aforementioned setting is satisfied.

 \subsection{Final part of the proof} Assume that there exists $x_0>0$ and $c<1$ such that for all $x\geq x_0$ we have $\lambda_x/\mu_x\leq 1+c/x$. Then, starting from some $k\geq k_0$, where $k_0$ is chosen to be sufficiently large, from \eqref{25} we have
 \[
 \begin{aligned}
 &\int_{k-1}^{k}P_{x+1}\mu_{x+1}\mathrm{d}x+\int_{k-1}^{k}P_{x-1}\mu_{x-1}\left(1+\frac{c_{x-1}}{x-1}\right)\mathrm{d}x\\
 &=\int_{k-1}^{k}P_x\mu_x\mathrm{d}x-\int_{k-1}^{k}P_x\mu_x\left(1+\frac{c_x}{x}\right)\mathrm{d}x,
 \end{aligned}
 \]
where $c_x\leq c<1$, from which it readily seen that we must obtain $\lambda_k^*/\mu_k^*\leq1+c^*/k$ for all $k\geq k_0$, where $c^*< c+(k_0-1)^{-1}<1$.

If we assume that there exists $x_0>0$ and $c>1$ such that for all $x\geq x_0$ we have $\lambda_x/\mu_x\geq 1+c/x$, then similarly to the above we arrive at $\lambda_k^*/\mu_k^*\geq1+c^*/k$ for all $k\geq k_0$, $c^*>1$.

Thus, we arrive at the reduction to the birth-and-death process. The rest of the proof fully repeats the final part of the proof, \cite[Section 3.3]{A}, with small changes adapted to the case.

Let $\mathcal{F}_\tau=\sigma\{Z(t), 0\leq t\leq\tau; X_i, Y_j, 1\leq i,j<\infty\}$.

If $\varphi(Z(\tau),\tau)$ is a positive process vanishing in $L^{1}$ as $\tau\to\infty$, then denoting by $\triangle Z_{\tau,\tau+\sigma}$ the increment of the process $Z(\tau)$ in the interval $[\tau, \tau+\sigma)$, for the process $Z^2(\tau)-\tau$ we have:
\[
\begin{aligned}
\mathsf{E}\left[Z^2(\tau+\sigma)-(\tau+\sigma)~|~\mathcal{F}_{\tau-}\right]=&\mathsf{E}\left[(Z(\tau-)+\triangle Z_{\tau, \tau+\sigma})^2-(\tau+\sigma)~|~\mathcal{F}_{\tau-}\right]\\
=&Z^2(\tau-)+2\mathsf{E}(Z(\tau-)\triangle Z_{\tau, \tau+\sigma}~|~\mathcal{F}_{\tau-})-\tau\\
&-\sigma+\mathsf{E}\left(\triangle Z_{\tau,\tau+\delta}^2~|~\mathcal{F}_{\tau-}\right)\\
=&Z^2(\tau)-\tau+o(Z(\tau)),
\end{aligned}
\]
since $\mathsf{E}(Z(\tau-)\triangle Z_{t,t+s}~|~\mathcal{F}_{\tau-})=o(Z(\tau))$, and
\[
\mathsf{E}(\triangle Z_{\tau,\tau+\sigma}^2|\mathcal{F}_{\tau-})=\mathsf{E}\triangle Z_{\tau,\tau+\sigma}^2=O(1)
\]
as $\tau\to\infty$. The last follows from Wald's identity \cite{F}. Specifically, we have $\mathsf{E}\triangle Z_{\tau,\tau+\sigma}^2\leq\mathsf{E}\sum_{i=1}^{N_\sigma}(X_i^2+Y_i^2)=\sigma(2+\mathsf{Var}(X_1)+\mathsf{Var}(Y_1))$, where $N_\sigma$ denotes the number of events in time $\sigma$ of a Poisson process with rate $1$. So, $\mathsf{E}Z^2(\tau)=\tau+o(\tau)$.

Along with this, we also have as follows: for any $\epsilon>0$ there exist $\delta>0$ and $\tau_0$ such that for all $\tau>\tau_0$
\[
\mathsf{P}\left\{\left|\frac{Z^2(\tau)}{\tau}-1\right|<\delta\right\}>1-\epsilon.
\]
Therefore, for the function $\varphi(Z(\tau), \tau)$ for large $\tau$, due to the continuity property of this function we have
$\varphi(Z(\tau), \tau)\asymp\varphi(Z(\tau), Z^2(\tau))$, where `$\asymp$' means that $\varphi(Z(\tau), \tau)/\varphi(Z(\tau), Z^2(\tau))$ converges to 1 in probability as $\tau\to\infty$.

Similarly, for the parameters $\lambda_{n,t}$ and $\mu_{n,t}$ considered in \cite{A} for large $\tau$ we have $\lambda_{Z(\tau),\tau}\asymp\lambda_{Z(\tau),Z^2(\tau)}$ and $\mu_{Z(\tau),\tau}\asymp\mu_{Z(\tau),Z^2(\tau)}$, which implies 
\[
\frac{\lambda_{Z(\tau),\tau}}{\mu_{Z(\tau),\tau}}\asymp\frac{\lambda_{Z(\tau),Z^2(\tau)}}{\mu_{Z(\tau),Z^2(\tau)}}.
\]
So, under the condition $Z(\tau)=n$ for large $n$, we obtain
\[
\frac{\lambda_{n,\tau}}{\mu_{n,\tau}}\asymp\frac{\lambda_{n,n^2}}{\mu_{n,n^2}},
\]
which is the required asymptotic presentation given in \cite{A}.

\section{Concluding remarks}\label{S6}

In this paper, we have developed the earlier result of \cite{A} for random walks taking the values in $\mathbb{R}$. The assumptions under which the study has been conducted are natural and less restrictive than those in \cite{MV}. Furthermore, the obtained result is more general than that in \cite{MV}. All the examples that have been considered in \cite{A} are automatically valid for the model considered in the present paper. In particular, we completely solved an open problem that has been formulated in \cite{MV} and partially solved in \cite{A} in terms of the present paper.

\subsection*{Acknowledgement} The author thanks Dr. Tuan-Minh Nguyen for useful conversation. As well, the author thanks all the people, who made critical comments at official or private level.


\subsection*{Disclosures and declarations}
No conflict of interest was reported by the author.


\end{document}